\documentclass[11pt]{article}
\usepackage[utf8x]{inputenc}
\usepackage{amsmath,amssymb,amsfonts,latexsym}

\usepackage{array}
\usepackage{booktabs}
\usepackage{multirow}
\usepackage{xcolor}

\usepackage{graphicx}   
\usepackage{pgfplots}
\pgfplotsset{compat=1.18}

\newtheorem{t1}{Theorem}[section]

\newtheorem{d1}{Definition}[section]

\begin{document}
	\title{Normalized Fractional Order Entropy-Based Decision-Making Models under Risk}
	\author{Poulami Paul\and Chanchal Kundu
		\and
		Department of Mathematical Sciences\\
		Rajiv Gandhi Institute of Petroleum Technology\\
		Jais 229 304, U.P., India}
	\date{January, 2026} 
	\maketitle
	\begin{abstract} 
		Constructing efficient portfolios requires balancing expected returns with risk through optimal stock selection, while accounting for investor preferences. In a recent work by Paul and Kundu (2026), the fractional-order entropy due to Ubriaco was introduced as an uncertainty measure to capture varying investor attitudes toward risk. Building on this foundation, we introduce a novel normalized fractional order entropy aligned with investors’ risk preferences that combines normalized fractional entropy with expected utility and variance. Risk sensitivity is modeled through the fractional parameter $ q $, interpolating between conservative or risk aversion ($ q\to 0 $) and adventurous or high risk tolerance ($ q\to 1 $) attitudes. Furthermore, the robustness and statistical significance of the fractional order entropy-based risk measure, termed normalized expected utility-fractional entropy (NEU-FE) and normalized expected utility-fractional entropy-variance (NEU-FEV) risk measures are explained with the help of machine learning tools, including Random forest, Ridge regression, Lasso Regression and artificial neural networks by using Indian stock market (NIFTY50). The results confirm that the proposed decision models support investors in making high-quality portfolio investments. 
	\end{abstract}
	\section{Introduction}
	The formalization of information theory by Shannon (1948) with the introduction of the entropy concept to quantify uncertainty of an action or event is considered one of the most significant contributions to research in statistics, computer science, electrical engineering and other interdisciplinary areas. Since the formulation of information-theoretic concepts, it has evolved largely for describing complex physical systems. One such important milestone was achieved by Yang and Qiu (2005) in decision theory through finding its use in modeling risky decisions by defining an expected utility-entropy (EU-E)-based risk measure. The EU-E model incorporated the subjective elements like decision makers’ attitudes towards risk or costs through an expected utility function and the objective uncertainty through the Shannon entropy function. The goal was to minimize the risk score measured in terms of the entropy and utility functions to identify the most stable, rational, or appropriate decision given the available information about the action space and the corresponding states of nature. Here, the risk is associated with the changes in the future value or cost of an asset or security, which an investor chooses to buy, due to fluctuations in market conditions or other uncertain events.
	
	\hspace*{0.2in} The EU-E risk measure is based on the general decision-making model which is laid upon the fact that people choose an action with low uncertainty and high expected utility. It served as a normative model of choosing profitable actions or logical decisions and sometimes also served as a descriptive model of economic behaviors. The general decision analysis model $\mathcal{G} = (\Theta, \mathcal{A},\vartheta)$ comprises of the state space represented by $\Theta = \{\theta\}$, the action space symbolized by $\mathcal{A}$ and the utility function $\vartheta(X)$ representing an individual risk preference or risk aversion attitude. The payoff function defined on the space $\mathcal{A} \times \Theta$ is denoted by $X = X(a,\theta).$ \\
	\hspace*{0.2in} Consequently, when the action space $\mathcal{A} = \{A_1,A_2,\dots,A_l\}$ and the state space $\theta_i$ corresponding to action $A_i$ represented by $\theta_i = \{\theta_{i1}, \theta_{l2},\dots,\theta_{lm_i} \}$ are considered to be finite, then we obtain a payoff $X = X(A_i,\theta_{ij}) = x_{ij}$ corresponding to action $A_i$ with state $\theta_{ij}$ with $\theta_i$ having a probability distribution $\{p_{ij} = P(X=x_{ij})\}$ which is the probability of occurrence of the payoff outcome $x_{ij}$, where $i = 1,2,\dots,l;\quad j = 1,2,\dots m_i.$ This combination of outcomes and their associated probabilities for the action $A_i$ can be written as (see Yang and Qiu, 2014 and references therein):
	\begin{equation}
		A_i = \begin{pmatrix}
			x_{i1} & x_{i2} & \dots & x_{im_i} \\
			p_{i1} & p_{i2} & \dots & p_{im_i}
		\end{pmatrix} 
	\end{equation}
	or can be expressed as $A_i = (x_{i1},	p_{i1};x_{i2},	p_{i2};\dots; x_{im_i},	p_{im_i})$. The tabulated form of the general decision model can be found in Yang and Qiu (2005). \\
	\hspace*{0.2in} However, with the rise in complexity and risk of actions or the attitude of individuals towards risk, the Yang and Qiu model failed to describe risky decisions properly. It was also observed that the EU-E model didnot agree with actual patterns of human behavior and indiviual risk perceptions (Fischer and Kleine, 2006). Further attempts were made to improve this measure and provide additional insights into decision problems such as the ones made by Marley and Luce (2008) and Luce et al. (2005). A significant contribution was made by Yang and Qiu (2014) to standardize the EU-E decision model by defining the normalized expected utility-entropy (NEU-E) based risk measure, where the regular Shannon entropy was replaced by normalized Shannon entropy. In the study by Yang and Qiu (2014), the Shannon entropy was normalized by dividing it with the maximum entropy value. The normalized Shannon entropy is defined as:
	\begin{equation}
		NS(p) = \begin{cases}
			-\sum_{x=1}^{n} p(x) \log p(x) /\log (n), &~\text{for}~ n  > 1 \\
			0 &~\text{for}~ n  = 1;
		\end{cases}
		\label{nsh}
	\end{equation} where $S(p) = -\sum_{x=1}^{n} p(x) \log p(x) $ is the Shannon entropy and $\log (n)$ is the  maximum value of Shannon entropy computed from actions with equally likely outcomes as in uniform distributions. Unless stated otherwise, $\log(\cdot)$ refers to the natural logarithm throughout this study. The normalization of entropy is particularly significant because it removes scale dependency, ensures comparability across different distributions, and enhances interpretability of uncertainty levels within the utility framework. This normalized measure (\ref{nsh}) provided reasonable explanations into solutions of risky decision problems such as the certainty effects which included two important empirical phenomena called the common ratio effect and the common consequence effect (Machina, 1987). The decision results observed in certainty effect described in Prospect theory by Kahneman and Tversky (1979) can also be interpreted meaningfully by NEU-E model. 
	Dong et al. (2016) further explored the expected utility and entropy-based decision-making model, highlighting the pivotal role of Shannon entropy in risk-based decision analysis. 
	Recently, Brito (2020) advanced the classical EU-E framework by integrating variance into the decision criterion, leading to the EU-EV model. The inclusion of variance played a crucial role in capturing the sensitivity of decisions to outcome dispersion, an aspect that entropy alone could not fully represent. This enhancement allowed the model to resolve several long-standing decision paradoxes, such as the Levy paradox (Levy, 1992) and the Allais paradox (Allais, 1953), which remained unexplained under both the EU-E and NEU-E formulations.
	Building upon this contribution, Brito (2022) further generalized the approach by incorporating normalized Shannon entropy (\ref{nsh}). This extension resulted in the normalized expected utility-entropy-variance (NEU-EV) model, offering a more comprehensive and balanced representation of investor behaviour by jointly accounting for utility, uncertainty, and variability in outcomes. Thus, this innovation aided in offering a more robust and interpretable framework for analyzing risk-uncertainty trade-offs. \\
	\hspace*{0.2in} Nevertheless, the existing decision models gives too much weightage on the utility function in the decision model to capture the perceived riskiness of decision-makers but it fails to explain which part of the individual risk attitude should be covered by the entropy. For instance, in Prospect Theory (Kahneman and Tversky, 1979), individual attitudes towards changes in wealth is covered by the value function while the ``pure attitude" towards risk is represented by the probability weighting function. Hence, the extant risk measures fall short in incorporating the dynamic nature of investors' risk attitudes in relation to the uncertainty of risky actions. They also fail to effectively utilize machine learning (ML) tools, for analyzing the non-linear relationships between the input and output response variables. Moreover, these models failed to explore the potentials of other generalized entropies of degree $q$, such as R\'enyi (R\'enyi, 1961), Tsallis (Tsallis, 1988) or fractional (Ubriaco, 2009) entropies, in supporting the decision-making models under risk influenced by varying risk attitudes of individuals. Besides, the selection of the utility function and the specification of the trade-off parameter in these models often lack empirical calibration or behavioral grounding. As a result, the relationship between an investor’s sensitivity to uncertainty and the expected utility derived from different choices is not fully justified or supported by real-world decision patterns. 
	More recently, Paul and Kundu (2026) proposed the fractional entropy-based measures of risk to address the gaps in the existing Shannon entropy-based decision models, termed the expected utility-fractional entropy (EU-FE) and expected utility-fractional entropy-variance (EU-FEV) risk measure. The decision models based on these measures successfully incorporated the subjective elements like an individual's sensitivity towards risk and uncertainty of a risky prospect represented by the expected utility and fractional entropy formulated by Ubriaco (2009). This study introduced the use of fractional entropy in developing decision-making models for risky and uncertain choices. It incorporated the risk attitudes of individuals into the decision-making framework represented by the fractional parameter $ q. $ The parameter $ q $ provides additional flexibility to the Ubriaco entropy to accommodate the unique risk attitudes of decision-makers towards the probability of occurrence of an outcome. This fractional entropy framework captures the behaviour of both risk-seeking and risk-averse individuals. Risk-seeking decision-makers, who exhibit an adventurous attitude, tend to assign greater importance to low-probability outcomes, reflecting their willingness to embrace uncertainty. In contrast, risk-averse individuals display a more conservative attitude, placing less weight on such unlikely events and expressing greater doubt about their occurrence (see Aggarwal, 2021a; b). 
	The FE-based risk measures are basically based on the general decision analysis models under risk, which includes different actions identified by their states of nature, with each state having its own characteristic distribution.
	
	\hspace*{0.2in} However, the fractional entropy proposed by Ubriaco can take any positive real value, which makes it difficult to compare uncertainty across actions or prospects that differ in the number or structure of their possible outcomes. In decision-making models, especially those involving multiple alternatives with varying outcome spaces. Such unbounded entropy values can distort the assessment of uncertainty and risk, leading to inconsistent or scale-dependent evaluations.
	To overcome this limitation, we normalize the Ubriaco fractional entropy so that its values lie within the fixed interval $[0,1]$. This normalization ensures that the entropy measure remains comparable across prospects regardless of their dimensionality, number of outcomes, or payoff structure. A bounded, unit-free entropy scale also enhances interpretability: values closer to $0$ represent low uncertainty, whereas values closer to $1$ indicate high uncertainty, enabling a clearer distinction among risky choices.
	By incorporating this normalized fractional-order entropy into the decision-making criterion, the present work introduces a new risk measure that unifies expected utility, variance, and fractional-order uncertainty in a consistent framework. This normalized formulation allows for more balanced risk assessment, avoids biases arising from unbounded entropy values, and ultimately provides a more robust tool for resolving decisions among risky prospects. 
	
	\hspace*{0.2in} The rest of the paper is organized into the following sections: Section 2 introduces the normalized fractional Ubriaco entropy. 
	Thereafter, we explore some fundamental properties of the normalized fractional entropy-based risk measures. In Section 3, we conduct a comparitive study of the performance of the proposed decision analysis models with the fractional entropy-based models through detailed statistical analysis with the help of ML tools such as Random Forest, Ridge Regression and Artificial Neural Network (ANN), highlighting their advantages over the existing fractional entropy models by employing real market data from Indian stock market index (NIFTY50). Finally, in Section 4, we conclude the present study.
	
	\section{Normalized fractional Ubriaco entropy}
	\begin{d1}
		If the distribution of a discrete random variable (r.v) is characterized by its probability mass function $ p(x) = P(X=x) $ for $ x^{th} $ outcome of an action or event with $n$ states of nature, then the fractional order entropy formulated by Ubriaco (2009) is defined as:
		\begin{equation}
			S_q (p) = \sum_{x=1}^{n}p(x) (-\log p(x))^q;~ q \in [0,1].
			\label{2.1}
		\end{equation} where $ s_x(p) = p(x) (-\log p(x))^q $ and $ I_q(p) = (-\log p(x))^q $ are treated as the fractional order entropy function and fractional order information gain function, respectively .
	\end{d1}
	Intuitively, when $ q = 1 $, the fractional order entropy (\ref{2.1}) is equivalent to the Shannon entropy expressed as: 
	\begin{equation}
		S(p) = -\sum_{x=1}^{n}p(x) \log p(x).
		\label{2.2}
	\end{equation}
	
	Ubriaco showed that the fractional entropy $ S_q(p) $ is always positive and the corresponding entropy function $ s_x(p) $ achieves a maximum value at $ p(x) = e^{-q}. $ Therefore, we get a loose upper bound for each term of the Ubriaco entropy given as $ s_x(p) = q^q e^{-q}. $ Thus, by summing up the maximum values of each $ s_x $ for $ n- $states of the system we obtain an upper bound of the entropy such that $ S_q \leq nq^q e^{-q}. $ Hence, the normalized fractional Ubriaco entropy can be defined as:
	\begin{equation}
		NS_q(p) = \frac{S_q(p)}{nq^q e^{-q}}
		\label{2.5}
	\end{equation} such that $ 0 \leq NS_q(p) \leq 1, ~ q \in [0,1]. $

\section{Normalized fractional entropy-based decision analysis models}
While evaluating risky alternatives, traditional expected utility entropy based risk measures often lead to inconsistencies, particularly when applied to classic paradoxes observed in empirical decision making experiments. A prominent example is the Allais paradox, which illustrates the certainty effect, which states the tendency of individuals to prefer a risk-free choice with sure gain having high probability over a riskier option, even when the latter offers a higher expected payoff. Other behavioural phenomena, such as the common ratio effect and the common consequence effect, also remain unexplained under conventional utility functions combined with Shannon entropy.
The normalized entropy-based risk  measures developed by Yang and Qiu (2014) and Brito (2022) successfully solved the inconsistencies in these decision problems to some extent. Here, the Shannon entropy was used to measure the objective uncertainty related to the risk or predictability of the occurrence of desired outcomes. However, the uncertainty is subjective to the risk tolerance levels or degree of conservatism of decision-makers. To take the subjective uncertainty into account, the fractional order entropy was incorporated into risk measures, where the fractional parameter was used to reflect the risk tolerance levels of individuals. The lower value of the fractional parameter $q$ close to zero was used to identify the conservative individuals, while the higher values closer to one represented the adventurous individuals (cf. Paul and Kundu, 2026).
The fractional order entropy-based decision model emphasized the potentials of fractional order entropy as an uncertainty measure in finance, both theoretically and empirically. It also investigated the strengths of fractional entropy in quantifying uncertainty in portfolio management by applying to the Portuguese stock market (PSI 20), collected from \textit{Investing.com}, verifying its effectiveness as an important component in risk measures and risk-based decision models. 

\hspace*{0.2in} Intrigued by the potentials of the fractional order entropy and motivated by the works of Yang and Qiu (2014) and Brito (2022), we intoduce a standardized form of the fractional entropy by normalizing it and integrate the proposed normalized entropy (\ref{2.5}) in the risk measure defined as follows:
\begin{d1}
	Consider a general decision analysis model $ G = (\Theta,\mathcal{A}_1,\vartheta) $, with an increasing utility function $ \vartheta = \vartheta(X(A_1,\theta)) $ corresponding to action $ A_1 \in \mathcal{A}. $ Then, the normalized expected utility- fractional entropy (NEU-FE) measure of risk associated with an action $ A_1 $ is defined by
	\begin{equation}
		R(A_1) = \begin{cases}
			\lambda NS_q(A_1,\theta)-(1-\lambda)\frac{E[\vartheta(X(A_1,\theta))]}{\max_{A_1\in \mathcal{A}}\{|E[\vartheta(X(A_1,\theta))]|\}} & \text{for}~ \max_{A_1\in A}\{|E[\vartheta(X(A_1,\theta))]|\} \neq 0\\
			\lambda NS_q(A_1,\theta)  & \text{for}~ \max_{A_1\in A}\{|E[\vartheta(X(A_1,\theta))]|\} = 0
		\end{cases}
		\label{3.8}
	\end{equation}
	where $ \lambda \in \mathcal{R} $ is the risk-tradeoff factor such that $ 0 \leq \lambda \leq 1 $ and $ NS_q(A_1,\theta) $ is the normalized fractional entropy of the distribution of the state of nature corresponding to action $ A_1. $ Here, we assume that $ max_{A_1\in A}\{|E[\vartheta(X(A_1,\theta))]|\} $ exists and is non-zero. 
	\label{def_3.1}	
\end{d1}
\begin{d1}
	Given a general decision analysis model $ G = (\Theta,\mathcal{A}_2,\vartheta) $, with an increasing utility function $ \vartheta = \vartheta(X(A_2,\theta)) $ corresponding to action $ A_2 \in \mathcal{A}_2. $ Then, the normalized expected utility- fractional entropy and variance (NEU-FEV) measure of risk associated with an action $ A_2 $ is defined by
	\begin{equation}
		R(A_2) = \frac{\lambda}{2}\Bigg[NS_q(A_2,\theta) + \frac{Var[X(A_2,\theta)]}{\max_{A_2\in \mathcal{A}_2}\{Var[X(A_2,\theta)]\}}\Bigg]-(1-\lambda)\frac{E[\vartheta(X(A_2,\theta))]}{\max_{A_2\in \mathcal{A}_2}\{|E[\vartheta(X(A_2,\theta))]|\}},
		\label{3.9}
	\end{equation}
	where $ \lambda \in \mathcal{R} $ is the risk-tradeoff factor such that $ 0 \leq \lambda \leq 1 $ and $ NS_q(A_2,\theta) $ is the normalized fractional entropy of the distribution of the state of nature corresponding to action $ A_2. $ Here, we assume that $ max_{A_2\in \mathcal{A}_2}\{|E[\vartheta(X(A_2,\theta))]|\} $ exists and is non-zero. When $ max_{A_2\in \mathcal{A}_2}\{|E[\vartheta(X(A_2,\theta))]|\} = 0, $ then $$ R(A_2) = \frac{\lambda}{2}\Bigg[NS_q(A_2,\theta) + \frac{Var[X(A_2,\theta)]}{max_{A_2\in \mathcal{A}_2}\{Var[X(A_2,\theta)]\}}\Bigg]. $$
	\label{def_3.2}
\end{d1}
\vspace{-1.5em}
\hspace*{0.2in} To assess decisions made by individuals across a finite set of actions based on a risk measure, the risk values ($ R(\mathcal{A}), \mathcal{A}= \mathcal{A}_1 ~\text{or}~ \mathcal{A}_2 $) associated with each action are compared. These risk values are determined by the subjective uncertainty of the outcomes and their perceived utility. The action with the lowest $ R(\mathcal{A}) $  value is selected. Therefore, the NEU-FE and NEU-FEV decision models can be redefined as:
\begin{d1}
	Let $ G = (\Theta,\mathcal{A},\vartheta) $ be a decision analysis model.
	\begin{enumerate}
		\item Let the action space $ \mathcal{A} = \mathcal{A}_1 ~\text{or}~ \mathcal{A}_2 $ consist of two actions $ A_{i1}, A_{i2} \in \mathcal{A}, ~i = 1,2 $ with related risk measures $ R(A_{i1}) $ and $ R(A_{i2}), $ respectively. Then the following results follow: \\ 
		$ (i) $ action $ A_{i1} $ is selected, i.e., $ A_{i1} > A_{i2} $, if $ R(A_{i1}) < R(A_{i2}); $ \\
		$ (ii) $ action $ A_{i2} $ is selected, i.e., $ A_{i1} < A_{i2} $, if $ R(A_{i1}) > R(A_{i2}); $ \\
		$ (iii) $ no action is superior i.e., any action can be selected or, $ A_{i1} = A_{i2} $, if $ R(A_{i1}) = R(A_{i2}). $ 
		\item If $ \mathcal{A} $ is a discrete set of $ n $ objects, i.e., $ \mathcal{A} = \{A_{i1},A_{i2},\dots,A_{in}\}, $ then the NEU-FE and the NEU-FEV order of preferences of actions will be in accordance with decreasing order of risk values. The most preferred action will then be 
		\[R(A_{ij}) = min_{A_{ij}\in \mathcal{A}} R(A_{ij}). \] 
		\label{def_3.3}
	\end{enumerate}
\end{d1} 
\vspace{-1.5em}
\hspace*{0.2in} Moreover, for $ \lambda = 0, $ the proposed decision models reduce to the expected utility model ($ A_{i1} > A_{i2} $ or $ R(A_{i1}) < R(A_{i2}) $ iff $E[\vartheta(A_{i1})] > E[\vartheta(A_{i2})] $). 
\subsection{Properties of the NEU-FE and NEU-FEV measures of risk}
The NEU-FE  and NEU-FEV risk measures obey the fundamental principles of risk measures, similar to those of previously established risk measures. This section will explore some of the key properties of these proposed risk measures.
\begin{t1}
	Let $ A_{11} $ and $ A_{12} $ be two different actions in the action space associated with a general decision model $ G = (\Theta,\mathcal{A}_1,\vartheta) $ having equal expected utility, $ NS_q(A_{11},\theta) $ and $ NS_q(A_{12},\theta)  $ denote normalized fractional entropies of the corresponding states of nature $\theta$, then $ R(A_{11}) < R(A_{12}) $ when $ NS_q(A_{11},\theta) < NS_q(A_{12},\theta) $ for some $ q. $ For the NEU-FEV model, $ R(A_{21}) < R(A_{22}) $ if and only if $$ \Big[NS_q(A_{21},\theta) + \frac{Var[X(A_{21},\theta)]}{max_{A_{21}\in \mathcal{A}_2}\{Var[X(A_{21},\theta)]\}}\Big] < \Big[NS_q(A_{22},\theta) + \frac{Var[X(A_{22},\theta)]}{max_{A_{22}\in \mathcal{A}_2}\{Var[X(A_{22},\theta)]\}}\Big]. $$ 
\end{t1}
\textbf{\textit{Proof.}}
The result immediately follows from Definitions \ref{def_3.1} and \ref{def_3.2}. $\hfill\square$\\

\hspace*{0.02in} Intuitively, according to the theorem mentioned above, an action with lower entropy and smaller variance will be associated with reduced risk.                       
\begin{t1}
	Let $ G = (\Theta,\mathcal{A},\vartheta) $ be a general decision analysis model, with the action space comprised of two actions $ A_{11} $ and $ A_{12}, ~ NS_q(A_{11},\theta) $ and $ NS_q(A_{12},\theta) $ denoting the normalized fractional entropies of the corresponding states of nature $ \theta $, respectively, then $ R(A_{11}) < R(A_{12}) $ when $ NS_q(A_{11},\theta) = NS_q(A_{12},\theta) $ for a particular value of $ q$ and $ E[\vartheta(X(A_{11},\theta))] > E[\vartheta(X(A_{12},\theta))]. $ Further, when $ NS_q(A_{21},\theta) = NS_q(A_{21},\theta) $, then $ R(A_{21}) < R(A_{22}) $ if and only if $E[\vartheta(X(A_{21},\theta))] > E[\vartheta(X(A_{22},\theta))]$. 	
	\label{thm_4.1}
\end{t1}
\textbf{\textit{Proof.}}
This result is also an immediate consequence of Definitions \ref{def_3.1} and \ref{def_3.2} and of the fact that for a decision model $ G(\Theta,\mathcal{A},\vartheta), $ we have $ Var(X(A,\theta)) = 0 $ if $  NS_q(A,\theta) = 0 $. $\hfill\square$\\

\hspace*{0.2in} From Theorem \ref{thm_4.1}, we get that the action having higher expected utility will be less riskier when the normalized fractional entropy values of the associated outcomes are equal. Further, the validity of the NEU-FE and NEU-FEV risk measures defined in (\ref{3.8}) and (\ref{3.9}) is examined by exploring several properties in relation to prior studies on decision models involving Shannon entropy and individual's risk perception.  

\begin{t1}
	Let $ G = (\Theta,\mathcal{A},\vartheta) $ be a general decision analysis model with non-negative outcomes:
	\begin{enumerate}
		\item Let the action space be described as $ \mathcal{A} = \{A, A+c\}, $ where $ c > 0 $ is a constant, then $$ R(A + c) < R(A), $$  
		i.e., risk declines when we add a positive constant to all outcomes of an action.
		\item If the action space contains two items, namely, $ \mathcal{A} = \{A,kA\}, $ where $ k > 1 $ is a constant, then for the NEU-FE decision model, $$ R(kA) < R(A), $$
		i.e., NEU-FE risk decreases if we multiply a positive constant to all outcomes of an action.\\
		\hspace*{0.2in} For the NEU-FEV decision model, $$ R(kA) < R(A), $$ if and only if \[\lambda \in \Bigg[0,~\frac{1-\frac{E[\vartheta(X)]}{E[\vartheta(kX)]}}{\frac{3}{2}-\frac{1}{2k^2}-\frac{E[\vartheta(X)]}{E[\vartheta(kX)]}}\Bigg).\]	\end{enumerate} 
\end{t1}
\textit{\textbf{Proof:}} \begin{enumerate}
	\item As we know that the utility function $ \vartheta(X) $ is increasing, then for $ c>0, $ we have $$ \vartheta(X+c) > \vartheta(X). $$ If $ \mathcal{A} = \{A,A+c\}, $ then $$ max_{A\in Y}\{|E[\vartheta(X(A,\theta))]|\} = E[\vartheta(X+c)] ~\text{and}~ Var(X+c)= Var(X). $$ Hence, from Definition \ref{def_3.2}, we get that 
	\[R(A_2) = \frac{\lambda}{2}\Big[NS_q(A_2,\theta) + 1\Big]-(1-\lambda)\frac{E[\vartheta(A_2)]}{E[\vartheta(A_2+c)]}\]
	and \[R(A_2+c) = \frac{\lambda}{2}\Big[NS_q(A_2+c,\theta) + 1\Big]-(1-\lambda).\]
	Similarly, from Definition \ref{def_3.1}, we have that
	\[R(A_1) = \lambda NS_q(A_1,\theta) - (1-\lambda)\frac{E[\vartheta(A_1)]}{E[\vartheta(A_1+c)]}\]
	and \[R(A_1+c) = \lambda NS_q(A_1+c,\theta)-(1-\lambda).\]
	Since, $ NS_q(A+c,\theta) = NS_q(A,\theta)$ ($ A = A_1 ~\text{or}~ A_2 $) for a fixed value of $ q \in [0,1], $ one gets that $ R(A+c) < R(A). $ 	
	\item As we have an increasing utility function and non-negative outcomes $ X = X(A_2,\theta) $ such that $ E[X] = 0, $ then we get $ \vartheta(kX) > \vartheta(X) $ for $ k>1. $ Since, $ \mathcal{A}_2 = \{A_2,kA_2\}, $ therefore
	\[max_{A_2\in \mathcal{A}_2}\{|E[\vartheta(X(A_2,\theta))]|\} = E[\vartheta(kX)] ~\text{and}~ Var(kX)= k^2E[X^2] > Var(X) = E[X^2].\]
	Hence, from Definition \ref{def_3.2}, we get that 
	\[R(A_2) = \frac{\lambda}{2}\Bigg[NS_q(A_2,\theta) + \frac{1}{k^2}\Bigg]-(1-\lambda)\frac{E[\vartheta(X)]}{E[\vartheta(kX)]}\]
	and \[R(kA_2) = \frac{\lambda}{2}\Big[NS_q(kA_2,\theta) + 1\Big]-(1-\lambda).\]
	Since entropy depends on the probabilities of the states of nature $ p(i) $ and $ q, $ for a particular value of $ q,~ NS_q(kA_2,\theta) = NS_q(A_2,\theta) $. Hence, it follows that $ R(kA_2) < R(A_2) $ if and only if
	\[\lambda\Bigg(\frac{3}{2}-\frac{1}{2k^2}-\frac{E[\vartheta(X)]}{E[\vartheta(kX)]}\Bigg)-\Bigg(1-\frac{E[\vartheta(X)]}{E[\vartheta(kX)]}\Bigg) > 0.\]
	Therefore, we have that for $ R(kA_2) < R(A_2) $ if and only if 
	\[\lambda \in \Bigg[0,~\frac{1-\frac{E[\vartheta(X)]}{E[\vartheta(kX)]}}{\frac{3}{2}-\frac{1}{2k^2}-\frac{E[\vartheta(X)]}{E[\vartheta(kX)]}}\Bigg).\]
	Similarly, it follows from Definition \ref{def_3.1} that
	\[R(A_1) = \lambda NS_q(A_1,\theta)-(1-\lambda)\frac{E[\vartheta(X)]}{E[\vartheta(kX)]}\]
	and \[R(kA_1) = \lambda NS_q(kA_1,\theta) -(1-\lambda).\]
	Thus, we obtain $ R(kA_1) < R(A_1) $ since for actions $ A_1 $ and $ kA_1, ~ NS_q(kA_1,\theta) = NS_q(A_1,\theta). $     
\end{enumerate}

\paragraph{Mathematical sensitivity of NEU-FEV to variance.}
Let \(r=(r_1,\dots,r_n)\) be the return sample for a given asset, with sample mean \(\bar r\) and sample variance
\[
s^2 \;=\; \frac{1}{n-1}\sum_{i=1}^n (r_i-\bar r)^2 .
\]
We consider the normalized-variance term \(\mathrm{NormVar}\) which (one common choice) can be the linear normalization
\(\mathrm{NormVar}=s^2 / V_{\max}\) or the z-standardized form \(\mathrm{NormVar}=(s^2-\mu_{s^2})/\sigma_{s^2}\). For simplicity of exposition below we use the linear normalization \(\mathrm{NormVar}=s^2/V_{\max}\); the qualitative conclusions below hold for other monotone normalizations as well.

The NEU–FEV target (for a fixed \(\lambda\in[0,1]\)) is
\[
T(r) \;=\; \frac{\lambda}{2}\big(\mathrm{NS}_q(r) + \mathrm{NormVar}(r)\big)\;-\;(1-\lambda)\,\mathrm{NEU}(r),
\]
where \(\mathrm{NS}_q\) is the normalized fractional entropy and \(\mathrm{NEU}\) the (normalized) expected utility. The partial sensitivity of \(T\) to the variance component is obtained by differentiation:
\[
\frac{\partial T}{\partial s^2}
\;=\; \frac{\lambda}{2}\,\frac{\partial \mathrm{NormVar}}{\partial s^2}
\;=\; \frac{\lambda}{2}\,\frac{1}{V_{\max}} \qquad\text{(for } \mathrm{NormVar}=s^2/V_{\max}\text{).}
\]
Thus, at the pointwise (local) level the contribution of variance to the target scales linearly with the prefactor \(\lambda/2\) and inversely with the normalization constant \(V_{\max}\). In words: any change \(\Delta s^2\) in the sample variance moves the target by
\[
\Delta T \approx \frac{\lambda}{2V_{\max}}\,\Delta s^2 .
\]

\medskip\noindent\textbf{Influence of a single extreme return.}
A more revealing view uses the influence of a contamination at a single observation \(r_i\). The influence function of the variance functional (population variance \(\sigma^2\)) is proportional to the centred squared deviation. For an observation \(x\) the influence is (up to scaling)
\[
\mathrm{IF}(x;\,\mathrm{Var}) \;=\; (x-\mu)^2 - \sigma^2.
\]
Hence, inserting an extreme value \(x\) (large \(|x-\mu|\)) changes the variance approximately by \((x-\mu)^2-\sigma^2\). The corresponding change in the target is
\[
\Delta T \approx \frac{\lambda}{2V_{\max}}\big[(x-\mu)^2-\sigma^2\big],
\]
which grows \emph{quadratically} with the extremeness \(|x-\mu|\). This quadratic dependence explains why outliers and heavy tails have an outsized effect on NEU-FEV but not on NEU-FE (which lacks the \(s^2\) term).

\medskip\noindent\textbf{Propagation into target variance.}
To quantify the effect of variance on the uncertainty of \(T\), compute the variance of \(T\) (treating \(\mathrm{NS}_q\), \(\mathrm{NormVar}\), and \(\mathrm{NEU}\) as random variables driven by the sample \(r\)):
\[
\operatorname{Var}[T]
= \left(\frac{\lambda}{2}\right)^2\operatorname{Var}[\mathrm{NS}_q]
+ \left(\frac{\lambda}{2}\right)^2\operatorname{Var}[\mathrm{NormVar}]
- 2\left(\frac{\lambda}{2}\right)(1-\lambda)\operatorname{Cov}[\mathrm{NormVar},\mathrm{NEU}]
+ \text{(other terms)},
\]
where the dominant new term introduced by augmentation is \(\left(\tfrac{\lambda}{2}\right)^2\operatorname{Var}[\mathrm{NormVar}]\). Using \(\mathrm{NormVar}=s^2/V_{\max}\), we have
\[
\operatorname{Var}[\mathrm{NormVar}] \;=\; \frac{1}{V_{\max}^2}\operatorname{Var}[s^2].
\]
For the sample variance the variance (asymptotically) depends on the fourth central moment \(\mu_4=\mathbb{E}[(r-\mu)^4]\). A standard expression (for independent draws) is
\[
\operatorname{Var}[s^2] \approx \frac{1}{n}\big(\mu_4 - \sigma^4\big) + o\!\left(\tfrac{1}{n}\right),
\]
so that
\[
\operatorname{Var}[\mathrm{NormVar}] \approx \frac{1}{nV_{\max}^2}\big(\mu_4 - \sigma^4\big).
\]
Consequently the variance contribution to \(\operatorname{Var}[T]\) scales as
\[
\left(\frac{\lambda}{2}\right)^2 \frac{1}{nV_{\max}^2}\big(\mu_4 - \sigma^4\big).
\]
Because \(\mu_4\) (the kurtosis-related moment) increases substantially for heavy-tailed distributions, the sampling variability of \(s^2\) (and hence of the target \(T\)) increases accordingly. This establishes that NEU-FEV inherits an extra term in its sampling variance that is proportional to the fourth moment of returns; NEU-FE lacks this extra source of uncertainty.

\medskip\noindent\textbf{Consequences (intuitive summary).}
\begin{itemize}
	\item The variance term enters the target linearly, but the sample variance itself depends quadratically on deviations and has sampling variance controlled by the fourth central moment. Therefore, NEU-FEV is \emph{more sensitive} to outliers and heavy tails than NEU-FE.
	\item An extreme return \(x\) perturbs NEU-FEV by an amount proportional to \((x-\mu)^2\), whereas NEU-FE is affected only via entropy or utility terms that typically scale more weakly with extreme deviations.
	\item In finite samples the extra contribution \(\propto(\lambda/2)^2\operatorname{Var}[s^2]\) increases the uncertainty of the target; as a result, models predicting NEU-FEV will generally exhibit larger prediction variance and (empirically) lower $R^2$ when compared to NEU-FE, unless the sample size \(n\) is large or the distribution is light-tailed.
\end{itemize}

\subsection{Statistical Validation of the proposed models using real data}
\begin{figure}[h]
	\centering
	\includegraphics[width=0.6\textwidth]{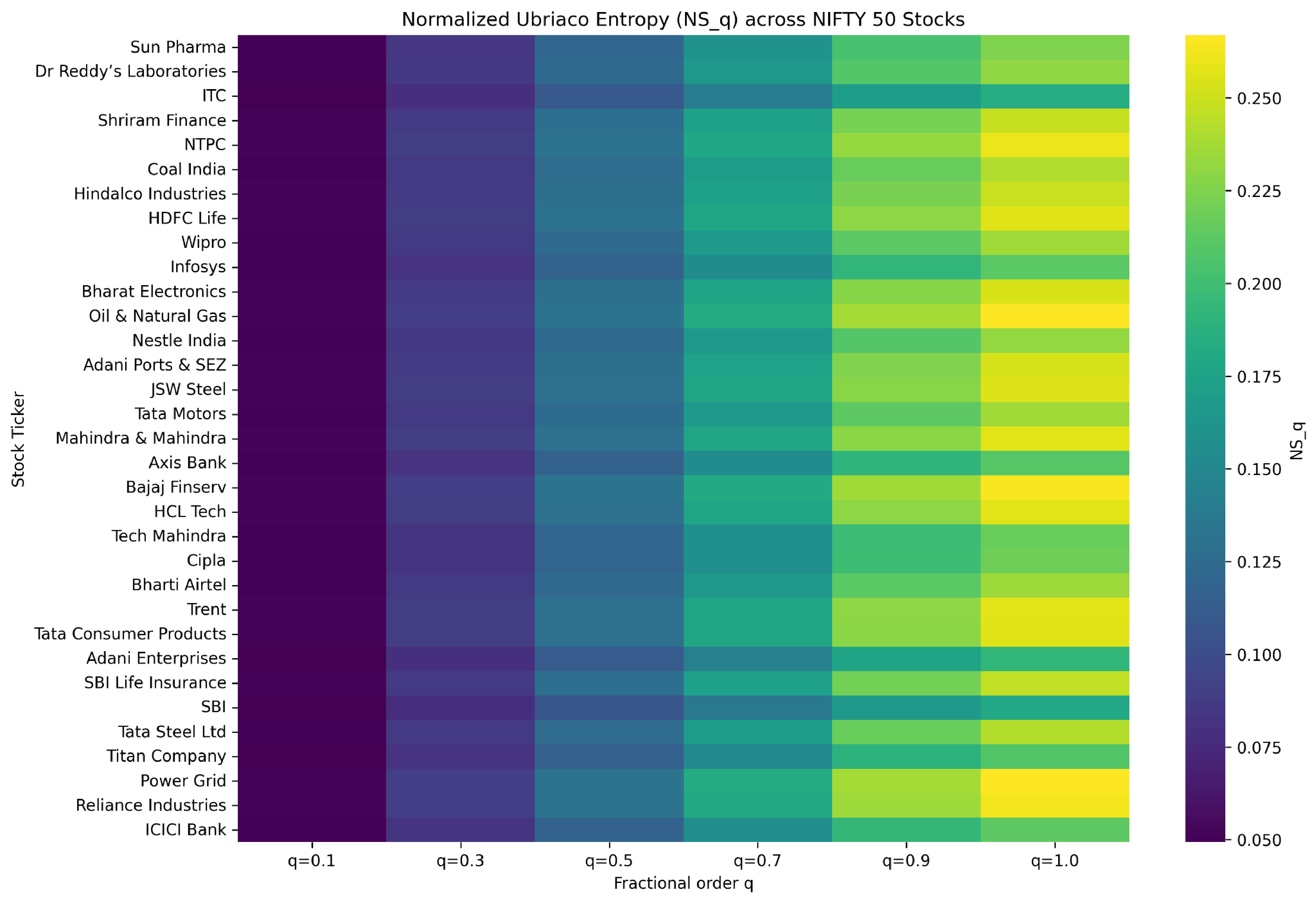} 
	\caption{Heatmap of variations of NEU-FEV measure with changes in $q$ for NIFTY50 stocks.}
	\label{FigH1} 
\end{figure}

\begin{figure}[h]
	\centering
	\includegraphics[width=0.5\textwidth]{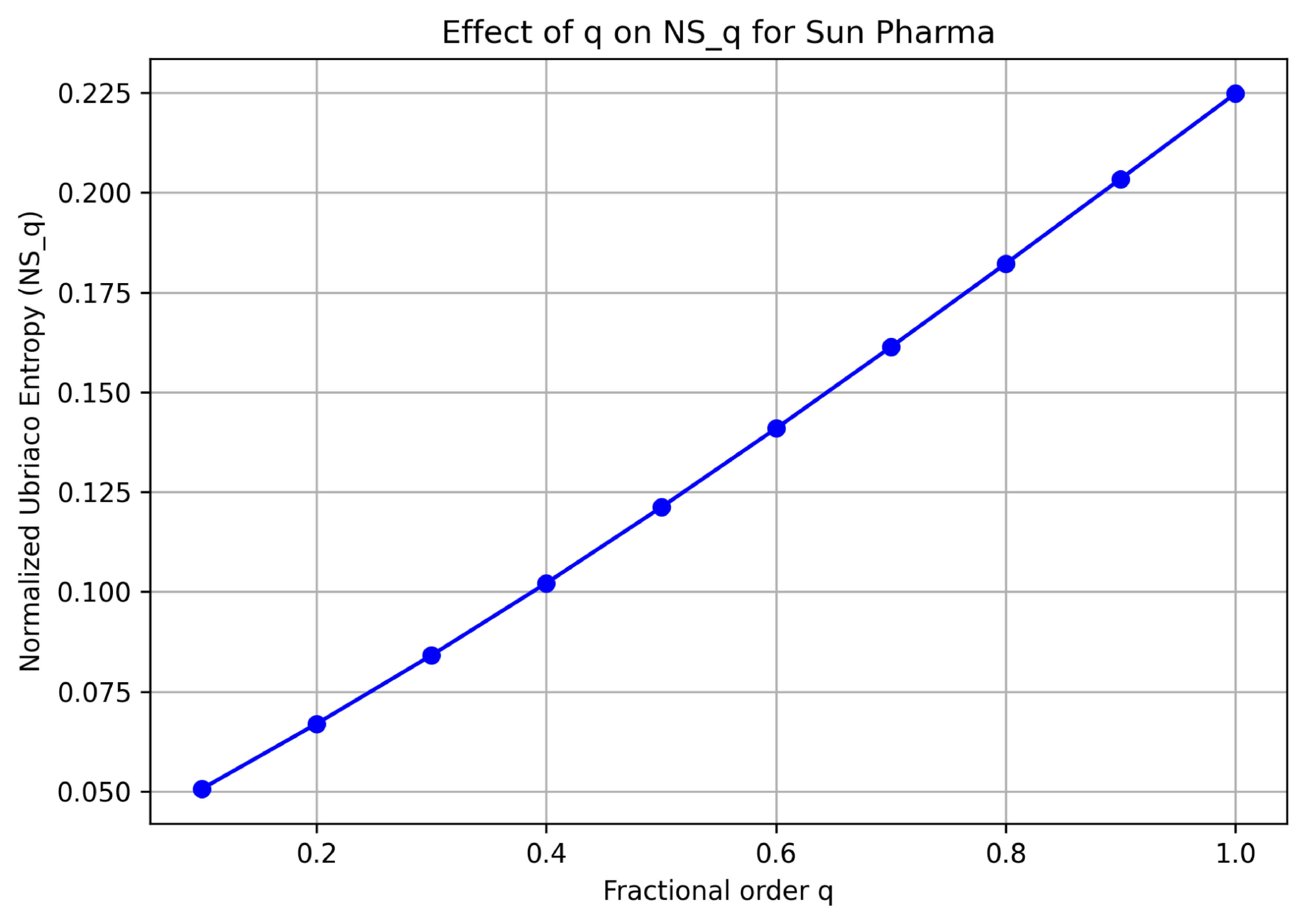} 
	\caption{Heatmap variations of NEU-FEV measure with changes in $q$ for a single NIFTY50 stock.}
	\label{FigH2} 
\end{figure}

The potentials of the fractional entropy-based decision analysis models in solving some important investment problems and decision paradoxes have recently been explored by Paul and Kundu (2026). The decision results reveal that it outperformed the existing Shannon entropy based decision models in most of the  cases and are consistent with the theoretical and empirical decision model results. In the present section, the strengths of the decision models based on the proposed normalized fractional entropy (\ref{2.5}) is explored and their superior performance in modeling decisions under risk is analyzed through some statistical tools using real data of  NIFTY 50 collected from \textit{Investing.com} by applying machine learning methods such as Random Forests, Ridge Regression and Artificial Neural Networks (ANNs). Furthermore, in order to analyze the behavior of the proposed models with respect to the model parameters $ q $ denoting the risk attitudes of decision-makers towards uncertainty, 
we conduct a detailed sensitivity analysis. The variations of the normalized fractional entropy-based risk measures with respect to $q$ are depicted in Figs.~\ref{FigH1} and \ref{FigH2}. As shown in these figures, the level of entropy and consequently the perceived uncertainty increases as the value of $q$ rises for each stock component. \\ 

\hspace*{0.2in} In the present study, we limit our comparative analysis to the most recent fractional entropy based risk measures, namely EU-FE and EU-FEV (cf. Paul and Kundu, 2026), as these models currently represent the strongest and most relevant alternatives within the fractional order decision making framework. The recent work has already established that fractional entropy based formulations consistently outperform the traditional Shannon entropy based EU-E and EU-EV criteria in resolving paradoxical choice patterns and capturing behavioural nuances in risky decision environments. Therefore in this section, to avoid redundancy and to ensure a meaningful and focused comparison, we benchmark our proposed normalized fractional entropy measures, NEU-FE and NEU-FEV, only against the existing fractional entropy-based counterparts, EU-FE and EU-FEV, as developed by Paul and Kundu (2026), using three predictive decision models: \textit{Random Forest}, \textit{Ridge Regression}, and \textit{Artificial Neural Networks (ANN)}. This shall serve as the appropriate state of the art baselines for evaluating further improvements.

\hspace*{0.2in} This study was performed using a one year stocks data of NIFTY 50 Index of the Indian market collected from \textit{Investing.com} from 04/03/2017 to 04/03/2018.  For this, we consider an investor willing to invest on $ l $ stocks components from $ s $ available choices. Then the action of selecting stock $ M_i $ will be defined on the action space $ \mathcal{A} = \{A_i\}, \quad i = 1,2,\dots,l $. The $l = 48$ available asset choices of the NIFTY 50 data are provided in Table \ref{Nifty}. We obtain a time series of 246 closing prices for the selected time interval corresponding to each stock $ S_i $ denoted by \{$ p_{i0},p_{i1},\dots,p_{iT} $\}, where $ T = 245. $ Then the returns for each stock $ S_i $ can be obtained as:
\begin{equation}
	r_{it} = \log \left(\frac{p_{it}}{p_{i(t-1)}}\right), \quad i = 1, \dots, 48; ~ t = 1,2,\dots,245.
	\label{4.9}
\end{equation} 
Then we get back a series of returns given by \{$ r_{i1},r_{i2},\dots,r_{i245}, $\} with respect to stock $ S_i $. Amongst all the available stock returns, we find that the minimum return value $ r_{min}= -0.187 $ and the maximum value $ r_{max}= 0.128 $ This gives us that 
$r_{it}\in [-0.187 ,0.128 ],$ for $ i = 1,2,\dots,48 $ and $ t = 1, \dots, 245. $ The interval [$ r_{min}, r_{max} $] is partitioned into $ \mathcal{J} = 15 $ sub-intervals $\mathcal{J}_k$, each with length $ \Delta = \frac{r_{max}-r_{min}}{\mathcal{J}} = 0.021 $. Therefore, we obtain 15 sub-intervals, listed as: $ \mathcal{J}_1 = [-0.187,-0.166), \mathcal{J}_2 = [-0.166,-0.145), \dots ,\mathcal{J}_{14} = [0.086,0.107], \mathcal{J}_{15} = [0.107,0.128].  $ 
Here, the probaility distribution or the relative frequency of the collected returns over $T$ previous days for the $i^{th}$  stock $S_i$ in the sub-interval $\mathcal{J}_k$ is computed as:
\begin{equation}
	p_{ik} = \frac{||r_{it} \in \mathcal{J}_k| t = 1,\dots, T||}{T}, ~i = 1,\dots, l, ~k = 1, \dots, \mathcal{J},
	\label{4.16}
\end{equation} where $||\cdot||$ represents the cardinality of a set and 
the $k^{th}$ sub-interval is obtained as:
\begin{equation}
	\mathcal{J}_k = \begin{cases}
		[r_{k-1},r_{k-1}+\Delta), & k = 1, \dots, \mathcal{J} \\
		[r_{k-1},r_k], & k = \mathcal{J}, 
	\end{cases}
	\label{4.17}
\end{equation} where we have $r_{max} = r_\mathcal{J}$. Then the expected return for the stock $S_i$ belonging to the sub-interval $\mathcal{J}$ collected over $T$ days, is evaluated as: 
\begin{equation}
	x_{ik} = \frac{1}{||r_{it} \in \mathcal{J}_k| t = 1,\dots, T||} \sum_{r_{it} \in \mathcal{J}_k,\\t = 1,\dots,T} r_{it}, ~ i = 1,\dots,l,~k =1,\dots,\mathcal{J}.
	\label{4.18}
\end{equation} Finally, by substituting the results from \eqref{4.16} and \eqref{4.17} into \eqref{4.9} and \eqref{4.18}, we obtain the risky action space represented by $A_i = (x_{i1},	p_{i1};\dots; x_{i15},	p_{i15};)$ of selecting the $i^{th}$ stock $S_i, i \in \{1,\dots,48\}$ from Indian market under NIFTY 50 index. The utility function used in this study is a popular S-shaped utility function, known to be concave for profits and convex for losses, as introduced by Kahneman and Tversky (1979) to represent the degree of risk aversion among decision-makers. The utility function is defined as:
\begin{equation}
	u(x) = \begin{cases}
		\log (1+x) & x \geq 0,\\
		-\log (1-x) & x < 0.
	\end{cases}
\end{equation} Consequently, the NEU-FE and NEU-FEV risk measures are defined using this utility function for each stock components. To check the performance of the model, all the stock components \{$ S_1,\dots,S_{48} $\} is chosen from the list of 48 stocks provided in Table \ref{Nifty}. The comparative results of NEU-FE and NEU-FEV risk measures with EU-FE and EU-FEV measures are summarized using mean squared error (MSE) and coefficient of determination ($R^2$) values in Table \ref{tab:decision_results}. Here, NEU-FE and NEU-FEV denote the proposed normalized measures, whereas EU-FE and EU-FEV serve as the existing baseline models.

\begin{table}[h!]
	\centering
	\caption{Stock components of NIFTY50}
	\begin{tabular}{|c|c|c|c|c|c|c|c|}
		\hline
		Stock & Components & Stock & Components & Stock & Components & Stock & Components \\
		\hline
		$S_1$ & ADEL & $S_{13}$ & REDY & $S_{25}$ & ITC & $S_{37}$ & SBI  \\
		$S_2$ & APSE & $S_{14}$ & EICH & $S_{26}$ & JSTL & $S_{38}$ & SHMF \\
		$S_3$ & APLH & $S_{15}$ & GRAS & $S_{27}$ & KTKM & $S_{39}$ & SUN\\
		$S_4$ & ASPN & $S_{16}$ & HCLT & $S_{28}$ & LART & $S_{40}$ & TCS \\
		$S_5$ & AXBK & $S_{17}$ & HDBK & $S_{29}$ & MAHM & $S_{41}$ & TACN \\ 
		$S_6$ & BAJA & $S_{18}$ & HDFL & $S_{30}$ & MRTI &$S_{42}$ & TAMO\\
		$S_7$ & BJFN & $S_{19}$ & HROM  & $S_{31}$ & NEST & $S_{43}$ & TISC \\
		$S_8$ & BJFS & $S_{20}$ & HALC & $S_{32}$ & NTPC & $S_{44}$ & TEML \\
		$S_9$ & BRTI & $S_{21}$ & HLL & $S_{33}$ & ONGC & $S_{45}$ & TITN  \\
		$S_{10}$ & BAJE & $S_{22}$ & ICBK & $S_{34}$ & PGRD & $S_{46}$ & TREN \\
		$S_{11}$ & CIPL & $S_{23}$ & INBK & $S_{35}$ & RELI & $S_{47}$ & ULTC\\
		$S_{12}$ & COAL & $S_{24}$ & INFY & $S_{36}$ & SBIL & $S_{48}$ & WIPR\\
		\hline
	\end{tabular}
	\label{Nifty}
\end{table}
\hspace*{0.2in} In this analysis, the stock market states or returns are first predicted using the designated machine learning techniques based on the chosen feature set. These predicted returns, together with the corresponding subjective utilities, are then used to compute the NEU-FE and NEU-FEV risk measures. The resulting decision outcomes illustrate how changes in $q$ influence the model's evaluation of risky prospects, providing insight into the role of fractional-order uncertainty in shaping investment decisions.

\subsection{Performance metrics comparison of the proposed Risk Measures}
This subsection presents the comparative evaluation of three supervised learning models—Lasso regression, Random Forest (RF), and Artificial Neural Network (ANN)—used to predict the normalized risk measure derived from the entropy–utility framework. Each model was trained on the feature matrix $X_{\mathrm{mat}}$ constructed from histogram-based entropy, generalized utility measures, variance-normalized dispersion statistics, skewness–kurtosis descriptors, and bootstrap-enhanced uncertainty indicators. The response variable $y$ represents the empirical risk metric associated with each financial asset.

\vspace{0.3cm}
\noindent\textbf{Influence of Regularization.}
The choice of regularization norm directly influenced both predictive behaviour and the decision-making model’s interpretability. The $L_1$ penalty in Lasso forced a data-driven feature selection that naturally identified the robust predictors of risk, thereby simplifying the resulting decision rules. In contrast, the $L_2$ penalty employed in Ridge (used here for comparison) stabilizes the coefficients without enforcing sparsity. Ridge tends to retain correlated entropy-derived features that Lasso suppresses, leading to different interpretations of feature importance and producing smoother risk gradients.

The combined findings show that Lasso enhances interpretability, RF emphasizes predictive strength and feature interactions, while the ANN offers the best flexibility in modelling the nonlinear geometry of entropy-based risk.

\begin{table}[htbp]
	\centering
	\caption{Comparisons of different risk measures. }
	\begin{tabular}{l | cc | cc}
		\toprule
		\multirow{2}{*}{\textbf{Model}} & \multicolumn{2}{c|}{\textbf{NEU-FE/EU-FE}}  & \multicolumn{2}{c}{\textbf{NEU-FEV / EU-FEV}} \\
		\cmidrule(lr){2-3} \cmidrule(lr){4-5} 
		& MSE & $R^{2}$ & MSE & $R^{2}$ \\
		\midrule
		Random Forest  & $\mathbf{8.44*10^{-5}}/0.0006$  & $\mathbf{0.9920}/0.9708$ & $\mathbf{0.0001} / 0.0002$ &  $\mathbf{0.9840}/ 0.9797$ \\
		Ridge   & $7.82*10^{-5}/7.41*10^{-5}$ & 0.9963/ 0.9966   & 0.0001 / 0.0001 & 0.9886/ 0.9915 \\
		Lasso  &  $2.54*10^{-6}/2.26*10^{-6}$ & $0.9999/0.9999$ & $\mathbf{1.12*10^{-5}}/1.19*10^{-5}$ & 0.9992/0.9994 \\
		ANN  & $\mathbf{0.0011}/0.0087$ & $\mathbf{0.8959}/0.6066$  &  $\mathbf{0.0011}/ 0.0099$  &  $\mathbf{0.8851}/ 0.3050$ \\
		\bottomrule
	\end{tabular}
	\label{tab:decision_results}
\end{table}


Table~\ref{tab:decision_results} presents a comparative performance analysis of four widely adopted machine-learning models such as Random Forest, Ridge regression, Lasso regression, and a two-layer Artificial Neural Network (ANN) for predicting the proposed fractional-entropy–based risk measures introduced in the next section. The evaluation is conducted using the Mean Squared Error (MSE) and the coefficient of determination ($R^{2}$), assessed for both the normalized entropy–utility–variance formulations (NEU-FE and NEU-FEV) and their non-normalized counterparts (EU-FE and EU-FEV). 

Across all risk formulations, Lasso regression consistently achieves the lowest MSE and the highest $R^{2}$ values, indicating near-perfect agreement between predicted and actual risk scores. For the NEU-FE and EU-FE models, Lasso attains MSE values on the order of $10^{-6}$ with $R^{2} = 0.9999$, demonstrating exceptional predictive precision and numerical stability. Its performance remains superior even for the variance-augmented formulations (NEU-FEV and EU-FEV), where the MSE stays within the $10^{-5}$ range while $R^{2}$ remains above $0.999$.

Ridge regression also delivers strong performance, achieving $R^{2} > 0.99$ across all fractional-entropy–based risk measures. Although its accuracy is slightly lower than that of Lasso, the Ridge results confirm that linear models with $L_{2}$ regularization effectively capture the structural behaviour of the entropy-driven decision models. The stability of Ridge regression further validates the smooth and well-conditioned nature of the NEU-FE and NEU-FEV formulations.

Random Forest demonstrates competitive performance, particularly within the normalized framework. For both NEU-FE and NEU-FEV, Random Forest achieves $R^{2}$ values in the range of $0.98$–$0.99$. However, its MSE values are larger than those of the regularized linear models, suggesting that the ensemble mechanism introduces additional variability when approximating the entropy-based risk landscape. This behaviour is consistent with the complex interaction patterns learned by tree-based ensembles, which may provide flexible modelling at the expense of noise amplification.

In contrast, the ANN model yields comparatively weaker results. Despite the capability of neural networks to approximate nonlinear functions, the results reveal substantially higher MSE values (between $10^{-3}$ and $10^{-2}$) and notably lower $R^{2}$ values, especially for the variance-extended risk measures. This suggests that the size and structure of the present dataset, together with the smooth analytical form of the risk measures, favour simpler parametric models over neural architectures. 

An important trend emerges when comparing normalized versus non-normalized risk formulations. For Random Forest, Ridge, and Lasso regression, the performance of the normalized variants (NEU-FE and NEU-FEV) is nearly identical to that of the EU-FE and EU-FEV models. This demonstrates that normalization preserves the essential structure of the decision model while improving numerical robustness. In several cases, particularly with Lasso and Ridge, the normalized measures show slightly higher $R^{2}$ values, underscoring the advantages of scale invariance and stable optimization introduced by normalization. The ANN exhibits larger discrepancies across normalized and non-normalized versions, further reinforcing that normalization primarily benefits linear models and tree-based ensembles.

Overall, Lasso regression emerges as the most reliable and accurate predictive model for both the normalized and non-normalized entropy-utility- variance decision frameworks. Its combination of sparsity, regularization-driven stability, and superior predictive fidelity makes it particularly well-suited for modelling fractional-entropy–based risk. Random Forest and Ridge regression follow as strong alternatives, while the ANN in its present configuration appears less suitable for representing the smooth, systematic structure inherent in the fractional-entropy–based risk surfaces.


\begin{figure}[h]
	\centering
	\includegraphics[width=0.6\textwidth]{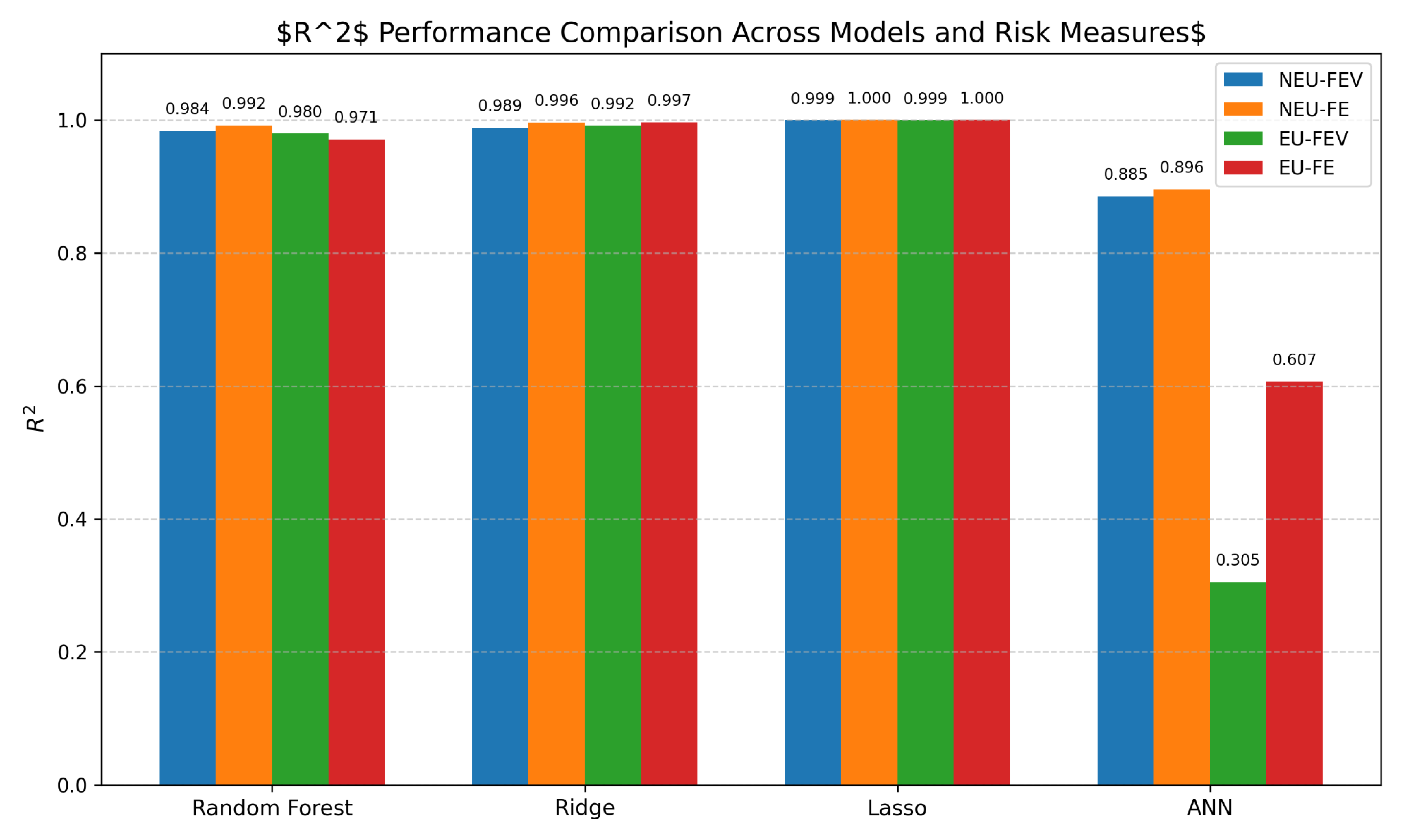} 
	\caption{Comparison of $R^2$ for all the fractional order entropy-based models.}
	\label{FigR2} 
\end{figure}

\section{Conclusion}
In this paper, we extend the fractional entropy-based decision-making framework of Paul and Kundu (2026) by introducing a normalized fractional-order entropy (NFE) and developing two associated risk-based decision models, NEU-FE and NEU-FEV. After establishing the key properties of these risk measures, we apply the proposed frameworks to real-world stock selection using NIFTY50 data by identifying the five most preferred stocks for a representative investor. 
A comprehensive comparative analysis is then performed against existing fractional entropy-based risk measures using widely adopted machine learning techniques, including Lasso and Ridge Regression, Random Forest, and a two-layer artificial neural network (ANN).

\hspace*{0.2in} Empirical results demonstrate that machine learning methods substantially enhance the modeling of fractional entropy-utility risk measures. Regularized linear models (Lasso and Ridge) offer strong predictive stability and interpretability by mitigating multicollinearity and highlighting the most influential entropy-utility components. Random Forest effectively captures nonlinear interactions and higher-order dependencies, while the ANN provides flexible function approximation but exhibits high error variability. Our findings further show that both normalized risk formulations, NEU-FE and NEU-FEV, retain or improve predictive performance relative to their non-normalized counterparts, confirming that normalization enhances numerical stability without altering the underlying decision structure. Moreover, incorporating variance increases model complexity without proportional predictive benefits, explaining the superior empirical performance of NEU-FE over NEU-FEV.

\subsection*{Data Availability Statement}
The details about the real market data consisting of 48 stock components of NIFTY50 are enlisted in Table \ref{Nifty} and can be collected from \textit{Investing.com}. 
\subsection*{Competing interests}
On behalf of all authors, the corresponding author declares that there is no conflict of interest or competing interests whatsoever in the preparation of the manuscript and conducting the study.


\begin{thebibliography}{}
	{\small
		\bibitem{am} Aggarwal, M. (2021a), Human decision-making through an entropic framework. \emph{Expert Systems with Applications}, \emph{183}. 
		\bibitem{a} Aggarwal, M. (2021b). Attitude-based entropy function and applications in decision-making. \emph{Engineering Applications of Artificial Intelligence}, \emph{104}, 104290. 
		\bibitem{al} Allais, M. (1953). The Behavior of the Rational Man in the Face of Risk, Critique of the Postulates and Axioms of the American School. \emph{Econometrica}, \emph{21}, 503--546. 
		\bibitem{ba} Brito, I. (2020). A decision model based on expected utility, entropy and variance. \emph{Applied Mathematics and Computation}, \emph{379}, 125285. 
		\bibitem{bb} Brito, I. (2022). The normalized expected utility – entropy and variance model for decisions under risk, \emph{International Journal of Approximate Reasoning}, \emph{148}, 174--201. 
		\bibitem{d} Dong, X., Lu, H., Xia, Y. and Xiong, Z. (2016). Decision-making model under risk assessment based on entropy. \emph{Entropy}, \emph{18}, 404.
		\bibitem{kt} Kahneman, D., and Tversky, A. (1979). Prospect Theory: An Analysis of Decision under Risk. \emph{Econometrica}, \emph{47}(2), 263--292. 
		\bibitem{fk} Fischer, K. and Kleine, A. (2007). Remarks on ``a measure of risk and a decision-making model based on expected utility and entropy" by Jipang Yang and Wanhua Qiu (EJOR 164 (2005), 792-799). \emph{European Journal of Operational Research}, \emph{182}, 469--474.
		\bibitem{le} Levy, H. (1992). Stochastic Dominance and Expected Utility: Survey and Analysis. \emph{Management Science}, \emph{38}(4), 555--593. 
		\bibitem{lu} Luce, R.D.,Ng., C.T., Marley, A.A.J. and Acz\'el (2005). Independence properties vis-a-vis several utility representations.. \emph{Risk and Decision Analysis}, \emph{58}, 77--143.
		\bibitem{ml} Marley, A.A.J. and Luce, R.D. (2008). Utility of gambling II: Risk, paradoxes, and data. \emph{Economic Theory}, \emph{36}, 165--187.
		\bibitem{m} Machina, M.J. (1987). Choice under uncertainty: Problems solved and unsolved. \emph{Journal of Economic Perspectives}, \emph{1}, 121--154.
		\bibitem{pk} Paul P. and Kundu, C. (2026). Fractional order entropy-based decision-making models under risk. \emph{Communications in Nonlinear Science and Numerical Simulation}, \emph{152}, 109106. 
		\bibitem{r} R$\acute{e}$nyi, A. (1961). On measures of entropy and information. \emph{In: Neyman, J., Ed., 4th Berkeley Symposium on Mathematical Statistics and Probability, Berkeley}, \emph{1}, 547-561.
		\bibitem{s} Shannon, C.E. (1948). A mathematical theory of communication. \emph{ The Bell System Technical Journal}, \emph{27}, 379--423, 623--656. 
		\bibitem{t} Tsallis, C. (1988). Possible generalization of Boltzmann-Gibbs statistics. \emph{Journal of Statistical Physics}, \emph{52}, 479--487. 
		\bibitem{u} Ubriaco, M.R. (2009). Entropies based on fractional calculus. \emph{Physics Letters A}, \emph{373}(30), 2516--2519. 
		\bibitem{yqa} Yang, J., and Qiu, W. (2005). A measure of risk and a decision-making model based on expected utility and entropy. \emph{European Journal of Operational Research}, \emph{164}(3), 792--799. 
		\bibitem{yqb} Yang, J., and Qiu, W. (2014). Normalized Expected Utility-Entropy Measure of Risk. \emph{Entropy}, \emph{16}(7), 3590--3604. 
	}
\end{thebibliography}
\end{document}